


\input amstex
\documentstyle{amsppt}

\baselineskip=12pt
\magnification=1200
\pagewidth{5.4in}
\pageheight{7.2in}
\parskip 12pt

\expandafter\redefine\csname logo\string@\endcsname{}
\NoBlackBoxes                
\NoRunningHeads
\redefine\no{\noindent}

\define\C{\Bbb C}
\define\R{\Bbb R}

\define\al{\alpha}
\define\be{\beta}

\define\ep{\varepsilon}
\define\Om{\Omega}
\define\De{\Delta}
\define\th{\theta}

\define\sub{\subseteq} 
    
\define\es{\emptyset}
     
\redefine\b{\partial}

\define\st{\ \vert\ }   

\redefine\ll{\lq\lq}
\redefine\rr{\rq\rq\ }
\define\rrr{\rq\rq}

\redefine\deg{\operatorname {deg}}
\define\Hol{\operatorname {Hol}}
\define\Map{\operatorname {Map}}

\define\SP{\operatorname {SP}}
\redefine\P{\operatorname {P}}
\redefine\F{\operatorname {F}}

\define\CC{\operatorname {C}}
\redefine\Q{\operatorname {Q}}
\define\jet{\operatorname {jet}}
\define\scan{\operatorname {scan}}
\define\Sec{\operatorname {Sec}}
\redefine\Re{\operatorname {Re}}

\define\CalSP{\Cal S\Cal P}

\topmatter
\title Spaces of polynomials with roots of bounded multiplicity\\
\endtitle
\footnote""{\it 1991 Mathematics Subject Classification. \rm 55P35, 58D15, 57R45}
\author M. A. Guest, A. Kozlowski, and K. Yamaguchi 
\endauthor
\abstract
We describe an 
alternative approach to some results of Vassiliev (\cite{Va1}) on
spaces of polynomials,  by using the
\ll scanning method\rr which was used by Segal (\cite{Se2}) 
in his investigation
of spaces of rational functions.  We
explain how these two approaches are
related by the Smale-Hirsch Principle or the h-Principle of Gromov.  We
obtain several generalizations, which
may be of interest in their own right.
\endabstract
\endtopmatter

\document

\head
\S 1 Introduction
\endhead

\noindent{\it Polynomials and rational functions}

The principal motivation for this paper derives from work of 
Vassiliev (\cite{Va1}, \cite{Va2}), in which he
describes a general method for calculating the cohomology of certain 
spaces of polynomial functions (and more generally,
\lq\lq complements of discriminants\rq\rq). 
As his paradigmatic
example, he takes the space $\P^d_n(\Bbb R)$  of real  
polynomials of the form 
$$
\spreadlines{2\jot}
x^d + a_{d-1}x^{d-1}+ \dots + a_0
$$
which have no $n$-fold real roots (but may  have complex ones
of arbitrary multiplicity!)  There is a \ll jet map\rr
$\P^d_n(\Bbb R)\to\Omega \Bbb R P^{n-1}$ given by
$f\mapsto [f;f^{\prime};\dots;f^{(n-1)}]$, and the image
of this map lies in a component $\Omega_{[d]} \Bbb R P^{n-1}$,
where $[d]=d$ mod $2$. One then has:

\proclaim{Theorem (Vassiliev)}
If $n\geq 4$, the jet map 
$\P^d_n(\Bbb R)\to\Omega_{[d]} \Bbb R P^{n-1}$ is a
homotopy equivalence up to dimension $([d/n]+1)(n-2)-1$.
\endproclaim

\no To say that a continuous map $f:X\to Y$ is a {\it homotopy
equivalence up to dimension} $d$ means that the induced homomorphism
$f_*:\pi_j(X)\to \pi_j(Y)$ of homotopy groups
is bijective when $j<d$ and surjective when $j=d$.

It follows from the theorem that 
the $d(n-2)$-skeleton of the space $\P^{nd}_n(\Bbb R)$ realises
the first $d+1$ cells in the well known cell decomposition
$\Omega_{[d]} \Bbb R P^{n-1} \simeq$
$\Om S^{n-1}\simeq $ {}
$e^0\cup e^{n-2}\cup e^{2(n-2)}\cup \dots$ of
the based loop space of the sphere. Homological considerations
show that $\P^{nd}_n(\Bbb R)$ has no homology above
dimension $d(n-2)$, so in fact  $\P^{nd}_n(\Bbb R)$ realises
these first $d+1$ cells exactly.

Vassiliev also considers the space $\SP^d_n(\Bbb C)$ of complex
polynomials without $n$-fold complex roots,  but states an analogue of
the above theorem only for
$n=2$ (and for homology groups), 
which gives a well known theorem of May and Segal 
(see \cite{Se1}) on the configuration space of {\it distinct}
points in $\C$. For arbitrary $n\ge 2$ one has
the following theorem, whose proof we shall give later on:

\proclaim{Theorem (Theorems 2.4 and 2.9)}
The jet  map 
$$
\SP^d_n(\Bbb C) \to \Omega_d^2\Bbb CP^{n-1},\quad
f\mapsto [f;f^{\prime};\dots;f^{(n-1)}]
$$
is a homotopy equivalence up to dimension $(2n-3)[d/n]$
if $n\ge 3$, and 
a homology equivalence up to dimension $(2n-3)[d/n]$
if $n=2$.
\endproclaim

\no In particular, for  $n\ge 3$, we have a homotopy equivalence
$$
\spreadlines{2\jot}
\lim_{d\to\infty}\SP^d_n(\Bbb C)
@>>> 
\Omega_0^2\Bbb CP^{n-1}.
\tag 1
$$ 

\remark{Remark}
We deliberately write $\SP^d_n(\Bbb C)$ --- rather 
than $\P^d_n(\Bbb C)$ --- in the complex case, because
$\SP^d_n(\Bbb C)$ is a subspace of the symmetric
product $\SP^d(\Bbb C)$ (the space of all complex 
polynomials of degree $d$). Note that  $\P^d_n(\Bbb R)$ is not
a subspace of $\SP^d(\Bbb R)$.
\endremark

To summarize, we may say that the space of real polynomials without
$n$-fold real roots is a model for $\Omega \R P^{n-1}$, and
the space of complex polynomials without
$n$-fold complex roots is a model for $\Omega^2 \C P^{n-1}$.
We shall show later that the space of complex polynomials without
$n$-fold real roots is a model for (the universal
covering space of) $\Omega \C P^{n-1}$,
and that the space of real polynomials without
$n$-fold complex roots is a model for the subspace of $\Omega^2 \C P^{n-1}$
consisting of maps $f:\C\cup\infty\to\C P^{n-1}$ such that 
$f(\overline z)=\overline{f(z)}$.

Let $\Hol_d(S^2,\Bbb CP^{n-1})$ 
denote the space consisting of all holomorphic (i.e. algebraic) maps 
$h:S^2\to \Bbb CP^{n-1}$ of degree $d$
which satisfy $h(\infty)=[1;0;\dots;0]$.
Concerning this space, one has the following theorem of \cite{Se2}:

\proclaim{Theorem (Segal)}
For $n\geq  2$ the inclusion $\Hol_d(S^2,\Bbb CP^{n-1}) 
\to  \Omega^2_d\Bbb CP^{n-1}$ is a homotopy
equivalence up to dimension  $d(2n-3)$.
\endproclaim

\no This theorem implies that 
$$ 
\spreadlines{2\jot}
\lim_{d\to\infty} 
\Hol_d(S^2,\Bbb CP^{n-1}) \simeq \Omega_0^2\Bbb CP^{n-1}.
\tag 2
$$
Combining (1) and (2), we obtain (for $n\ge 3$) the existence of a
homotopy equivalence
$$
\spreadlines{2\jot}
\lim_{d\to\infty}
\SP^d_n(\Bbb C) \simeq 
\lim_{d\to\infty}
\Hol_d(S^2,\Bbb CP^{n-1}).
\tag 3
$$ 
It turns out that there is an explicit description of this 
homotopy equivalence:

\proclaim{Proposition (Corollary 2.8)}
The jet map 
$$
\spreadlines{2\jot}
\SP^d_n(\Bbb C) \to \Hol_d(S^2,\Bbb CP^{n-1}),\quad
f\mapsto [f; f^{\prime}; \dots;f^{(n-1)}]
$$ 
induces the homotopy equivalence (3). 
\endproclaim

\no It was proved by Vassiliev (\cite{Va1}, \cite{Va2}) that,
if $n\geq 2$, there is a stable homotopy equivalence between
$\SP_n^{nd}(\Bbb C)$ and $\Hol_d(S^2,\Bbb CP^{n-1})$.
The case $n=2$ of this result was first proved in \cite{CCMM}.  
Proofs of the general case have been given subsequently in 
\cite{GKY2} and \cite{Kl2}.

Broadly speaking, two methods of proving theorems of the
above types appear in the literature.  One may be described
as \ll comparison of spectral sequences\rrr, after defining a
suitable filtration of each of the spaces concerned.
This was introduced in \cite{Va1}, \cite{Va2} and 
independently in \cite{CCMM}.  The
other uses a \ll scanning construction\rrr; this method,
due to Segal,
was developed in \cite{Mc1}, \cite{Mc2}, \cite{Se2}, \cite{Bo},
\cite{GKY1}, and \cite{Gu1}.
We shall use the scanning method here to give proofs of
the above results (including Vassiliev's original theorem).
We shall also prove various further generalizations, which we describe next.

\noindent{\it Equivariant homotopy equivalences}

The spaces $\SP^d_n(\Bbb C)$, $\P^d_n(\Bbb R)$ are examples of
a more general construction. Let $X,Y$ be subspaces of $\C$, and
let $\P^d_{Y,n}(X)$ denote the space of complex
monic polynomials $f$ of degree $d$ such that (i) $f(X)\sub X$
and (ii) $f$ has no $n$-fold roots in $Y$. Thus, 
$\P^d_{\C,n}(\C)=\SP^d_n(\Bbb C)$ and $\P^d_{\R,n}(\R)=\P^d_n(\Bbb R)$,
and we have seen that these spaces provide finite-dimensional
models for $\Omega^2\Bbb CP^{n-1}$ and
$\Omega\Bbb RP^{n-1}$ respectively. 

In a similar way, we
shall see that the spaces $\P^d_{\C,n}(\R)$ and $\P^d_{\R,n}(\C)$
provide models for loop spaces (as mentioned above). 
Now, $\R$ is the fixed point set
of the involution $\th:\C\to \C$, $\th(z)=\bar z$, and this involution
extends naturally to $\P^d_{\C,n}(\C)$ and $\P^d_{\R,n}(\C)$, and
to the corresponding loop spaces.  Using $\th$, the four results concerning
$\P^d_{Y,n}(X)$ (with $X,Y=\R$ or $\C$) can be summarized as
follows:

\proclaim{Theorem (Theorem 3.7)}
There is a $\th$-equivariant homotopy equivalence 
$$
\lim_{d\to\infty}
\P^d_{Y,n}(\C) \to \Map(Y^+,S^{2n-1})
$$
(if $n\ge 4$), where
$Y=\R$ or $\C$,  $Y^+=\R\cup\{\infty\}$ or $\C\cup\{\infty\}$, and where
$\Map$ indicates
the space of basepoint preserving continuous maps.
\endproclaim

By \cite{JS}, a $G$-map $\phi:A\to B$ of $G$-spaces is an equivariant
homotopy equivalence if and only if the maps
$\phi^{H}:A^{H}\to B^{H}$ of fixed point sets are
homotopy equivalences for all subgroups $H$ of $G$. In the
above theorem, $G$ is a group with two elements, so it has
precisely two subgroups.

In \cite{Se2}, Segal pointed out that the map 
$\lim_{d\to\infty} 
\Hol_d(S^2,\Bbb CP^{n-1})\to \Omega_0^2\Bbb CP^{n-1}$
is an equivariant homotopy equivalence in the same
way; this is equivalent to two statements, one for
complex rational functions and one for real rational
functions.  In this case also we have two additional \ll mixed\rr
spaces (which were not discussed by Segal).

\noindent{\it Holomorphic maps}

Another generalization of Vassiliev's theorem
is obtained by imposing \ll conditions of
bounded multiplicity\rr on the polynomials appearing in
Segal's theorem.  If we change the basepoint condition
on $f\in\Hol_d(S^2,\C P^{n-1})$ to $f(\infty)=[1;\dots;1]$,
then $f$ corresponds to an $n$-tuple $(p_1,\dots,p_n)$
of complex monic polynomials of degree $d$, such that
$p_1,\dots,p_n$ have no common root.  For $m\ge 2$, we may
define a subspace $\Hol^m_d(S^2,\C P^{n-1})$ of 
$\Hol_d(S^2,\C P^{n-1})$ by imposing
the additional condition that each $p_i$ belong to $\SP^d_m(\C)$.
The scanning construction then shows that $\Hol^m_d(S^2,\C P^{n-1})$
is a model for the double loop space $\Om^2 A_{n,m}$, where
$$
\spreadlines{2\jot}
A_{n,m}=\{(v_1,\dots,v_n)\in (\C^m-\{0\})^n
\st ((v_1)_1,\dots,(v_n)_1)\ne (0,\dots,0)\}.
$$
The space $A_{n,m}$ is an example of (the complement of) a
\ll subspace arrangement\rrr;  the topology of such spaces
has been studied intensively (see, for example, part III of \cite{GM},
as well as \cite{Va1}).

Segal's method extends to the case of
$\Hol(S^2,X)$, where $X$ is a toric variety (see
\cite{GKY1}, \cite{Gu1}, \cite{Gu2}). By imposing conditions of bounded
multiplicity on these spaces, we obtain many further results of
the above type. In each case, a subspace of $\Hol(S^2,X)$ gives
a (finite-dimensional) topological approximation 
to the double loop space of the
complement of a certain subspace arrangement.  The
\ll equivariant\rr results mentioned above also extend to
these examples.

The above results may be generalized in a different way
by replacing $\C$ by an open Riemann surface $\Sigma$, and
$S^2(=\C^+)$ by the one-point compactification $\Sigma^+$. One
obtains (for example) an equivalence up to some dimension
between $\SP^d_n(\Sigma)$ and the space $\Map(\Sigma^+,\C P^{n-1})$
of based maps. (Here, $\SP^d_n(\Sigma)$ is interpreted as a
subspace of the $d$-th symmetric product $\SP^d(\Sigma)$.)
For simplicity we shall restrict our exposition in this paper to the case
$\Sigma=\C$, referring to \cite{Mc1} for the method of
extension to general $\Sigma$ (and to open manifolds $\Sigma$
of arbitrary dimension).  However, it should be emphasized that
the extension to
manifolds other than $\C$ appears to be an advantage of the scanning
method.

Finally, we mention that the space $\Hol(S^2,X)$ may be identified
set-theoretically (although not topologically) with the space 
\ll $\Hol(S^1,X)$\rr which plays a role in the Gromov-Floer theory
of holomorphic curves (see for example \cite{Fu}).  It is noted
in section 3 of \cite{CJS} that 
$\Hol(S^1,\C P^{n-1})$ is homotopy equivalent to the universal
covering space of $\Om\C P^{n-1}$.  This fact can be proved 
by the methods described here, using $n$-tuples $(p_1,\dots,p_n)$
of polynomials which have no common root on $S^1$.

\noindent{\it The h-Principle}

The \ll scanning method\rr applies to all these examples because
in each case the space of polynomials involved can be identified
with a certain space of \ll labelled configurations\rrr.
By a labelled configuration we mean a finite set of distinct
points (usually in $\C$) where each point $z$ is labelled by an element
$m$ of a fixed partial monoid $M$.
The set of all such labelled configurations is topologized in
the usual way, except that two labelled points $(z_1,m_1)$,
$(z_2,m_2)$ are allowed to move towards each other and coalesce
\footnote{This is somewhat different to the usual notion of
labelled configuration in topology, where $M$ is simply a set.
In that case, distinct labelled points are never allowed to
coalesce.},
producing a
new labelled point $(z,m_1+m_2)$, if the sum $m_1+m_2$ exists in $M$.
More generally still, Kallel (\cite{Kl1}) has formulated the
notion of \ll particle space\rrr.

A deeper explanation for the existence of results of the above
type is suggested by Vassiliev, 
in terms of the {\it Smale-Hirsch Principle.} In 
its most general form, this is also known as the 
{\it h-Principle} of Gromov (see \cite{Gr}).
The relevant version of this principle says 
that --- under certain conditions --- 
the space of maps $M\to N$ whose
$k$-jets  avoid a \ll discriminant
variety\rr $S$ in the jet space $J^k(M,N)$ is homotopy
equivalent to the space consisting of all sections of the
bundle $J^k(M,N)-S\to M$. 
We shall describe a precise relation between the scanning method
and the h-Principle. The existence of such
a relation should not be surprising, as
the scanning method was in fact based on earlier ideas of Gromov (see
\cite{Mc2}).  However, the connection with
Gromov's work has been obscured in recent years by an emphasis 
(in the algebraic 
topology literature) on configuration spaces, so it seemed
worthwhile to explain here the original point of view. Indeed,
it might be argued that the h-Principle gives
the most natural approach to all \ll stable\rr  results of the
type considered here.

The paper is organized as follows. In \S 2 we discuss for
simplicity only the fundamental examples $\SP^d_n(\Bbb C)$ and
$\Hol_d(S^2,\Bbb CP^{n-1})$, and the relation between them. 
Various modifications of these examples (including Vassiliev's
original situation) are described in \S 3. Finally, in 
\S 4, we present the most general situation, and we explain
the relation with the h-Principle.

\no{\it Acknowledgements:}  We thank Sadok Kallel for sending us his preprints 
\cite{Kl1}, \cite{Kl2}, and Dai Tamaki for informing us about
the work \cite{Kt} of Fumiko Kato. Similar results to our
Theorem 2.2 were obtained independently both in \cite{Kl2} and in \cite{Kt}.
The first author was partially supported by a grant from the US National Science
Foundation, and the third author by a grant from the Ministry of
Education of Japan.

$${}$$

\head
\S 2 The fundamental example
\endhead

\noindent{\it Basic definitions}

For any space $X$, we denote by $\SP^d(X)$ the $d$-th symmetric
product of $X$.  By definition, this is the quotient space
$X^d/\Sigma_d$, where the symmetric group $\Sigma_d$ acts on
the $d$-fold product $X^d$ in the natural way.  An element of
$\SP^d(X)$ may be identified with a formal linear combination
$\al=\sum_{i=1}^k d_ix_i$, where $x_1,\dots,x_k$ are distinct
points of $X$ and $d_1,\dots,d_k$ are positive integers
such that $\sum_{i=1}^k d_i=d$.  We shall refer to $\al$ as
a \ll configuration\rr of points, the point $x_i$ having
multiplicity $d_i$.

In this section we shall be concerned with a subspace $\SP^d_n(X)$
of $\SP^d(X)$, defined as follows:

\proclaim{Definition 2.1}  For $n\ge 2$,
$\SP^d_n(X)=\{\sum_{i=1}^k d_ix_i
\in\SP^d(X) \ \vert\ d_i<n\ \text{for all}\ i\}$.
\endproclaim

\noindent Thus, $\SP^d_n(X)$ is obtained by imposing a
condition of \ll bounded multiplicity\rrr, namely that
all points $x_i$ (of any configuration) have multiplicity less than $n$. 
There is a filtration
$$
\spreadlines{2\jot}
\CC_d(X)=\SP^d_2(X)\sub \SP^d_3(X)\sub \dots \sub \SP^d_{d+1}(X)=\SP^d(X)
$$
where $\CC_d(X)$ denotes the space of \ll configurations of $d$
distinct points\rr in $X$.

If $A$ is a closed subspace of $X$, we define
$$\spreadlines{2\jot}
\SP^d_n(X,A)=
\{\sum_{i=1}^k d_ix_i\in\SP^d_n(X)\ \vert\ 
d_i<n\ \text{if}\ x_i\in X-A\}/\sim
$$
where $\al\sim\be$ if and only if $\al\cap (X-A)=\be\cap (X-A)$.
Thus, for $\SP^d_n(X,A)$, points in $A$ are \ll ignored\rrr.
When $A$ is nonempty, there is a natural inclusion map
$$\spreadlines{2\jot}
\SP^d_n(X,A)\to \SP^{d+1}_n(X,A)
$$
given by \ll adding a fixed point in $A$\rrr.  We define
$$\spreadlines{2\jot}
\SP_n(X,A)=\bigcup_{d\ge1}\SP^d_n(X,A).
$$ 
As a set, $\SP_n(X,A)$ is bijectively equivalent to the disjoint
union $\cup_{d\ge0}\SP^d_n(X-A)$,
but these two spaces are not in general homeomorphic.  For example, if
$\sum_{i=1}^k d_ix_i\in \SP_n(X,A)$, with 
$x_1,\dots,x_k\in X-A$, then a point $x_i$ can move into $A$
and \ll disappear\rrr.  In particular, $\SP_n(X,A)$ is connected,
if $X$ is connected.

We shall usually take $X$ to be an open subset of the
complex numbers $\C$.
Note that $\SP^d_n(\C)$ can be identified with the space of
complex polynomials of degree $d$ which are {\it monic},
and all of whose roots have {\it multiplicity less than} $n$.  (The
polynomial $\prod_{i=1}^k(z-x_i)^{d_i}$ corresponds to
$\sum_{i=1}^k d_ix_i$.)  In this special situation, 
there is a \ll stabilization map\rr
$$\spreadlines{2\jot}
\SP^d_n(\C)\to \SP^{d+1}_n(\C)
$$
which may be defined (up to homotopy) as follows.
If $U_d=\{z\in\C\ \vert\ \vert z\vert<d\}$, it is obvious
that $\SP^d_n(\C)$ is homeomorphic to $\SP^d_n(U_d)$.
Via this identification, the required map
$$\spreadlines{2\jot}
\SP^d_n(U_d)\to \SP^{d+1}_n(U_{d+1})
$$
is defined by $\sum_{i=1}^k d_ix_i\mapsto x+\sum_{i=1}^k d_ix_i$,
where $x$ is a fixed point of $U_{d+1}-U_d$.

\noindent{\it The scanning construction for configuration spaces}

To investigate the space $\lim_{d\to\infty}\SP^d_n(\C)$, we
use the \ll scanning map\rr
$$\spreadlines{2\jot}
s^d_n:\SP^d_n(\C)\to \Omega^2 \SP_n(\bar U,\partial \bar U)
$$
where $U=\{z\in\C\ \vert\ \vert z\vert < 1\}$. To define
this, we write 
$U(w)=\{z\in\C\ \vert\ \vert z-w\vert < \ep\}$,
where $\ep>0$ is fixed.
Let $\al=\sum_{i=1}^k d_ix_i\in\SP^d_n(\C)$. Then the
map 
$$\spreadlines{2\jot}
s^d_n(\al):\C\cup\infty\to \SP_n(\bar U,\partial \bar U)
$$
is defined by
$$\spreadlines{2\jot}
z\mapsto \al\cap\bar U(z) \in \SP_n(\bar U(z),\partial \bar U(z))
\cong \SP_n(\bar U,\partial \bar U).
$$
Note that $s^d_n(\al)$ is a basepoint-preserving map: the
point $\infty$ is always mapped to the empty configuration
in $\SP_n(\bar U,\partial \bar U)$.

As $\SP^d_n(\C)$ is connected, the image of $s^d_n$ lies
in a connected component of $\Omega^2 \SP_n(\bar U,\partial \bar U)$,
which we denote by $\Omega^2_d \SP_n(\bar U,\partial \bar U)$.
We shall see later
\footnote{This fact, as well as our Theorem 2.2, was noted
independently by Kallel (\cite{Kl2}) and Kato (\cite{Kt}).}
that 
$\SP_n(\bar U,\partial \bar U)\simeq\C P^{n-1}$; it is then easy
to show that $\Omega^2_d \SP_n(\bar U,\partial \bar U)$ is the
\ll$d$-th component\rr in the usual sense.  The main
reason for introducing the map $s^d_n$ is:

\proclaim{Theorem 2.2}  For $n\ge 3$, 
$$\spreadlines{2\jot}
\lim_{d\to\infty} s^d_n:  \lim_{d\to\infty} \SP^d_n(\C)
\to   \lim_{d\to\infty}\Omega^2_d\SP_n(\bar U,\partial \bar U)
\simeq \Omega^2_0\SP_n(\bar U,\partial \bar U)
$$
is a homotopy equivalence. For $n=2$, $\lim_{d\to\infty} s^d_n$
is a homology equivalence.
\endproclaim

\demo{Proof}  The proof is similar to the argument of \S 3 of
\cite{Se2}.  A detailed exposition,
more suited to the purposes of the present article, is given in 
Proposition 3.1 of \cite{Gu1}. Since this is an important argument,
however,
which will reappear in \S 4 in connection with the h-Principle,
we shall sketch the main ideas here (cf. \cite{Mc2}).

For real numbers $x,y\ge0$, let 
$D_{[x,y]}=\{z\in\C\ \vert\ x\le\vert z\vert\le y\}$. Consider the 
commutative diagram
$$
\CD
\SP_n(D_{[0,1]},D_{[\frac23,1]}) @>r>>  
\SP_n(D_{[0,1]},D_{[0,\frac13]}\cup 
D_{[\frac23,1]})\\
@Vs_1VV @VVs_2V\\
\Map(D_{[0,\frac23]},\SP_n(\bar U,\b\bar U)) @>r^\prime>> 
\Map(D_{[\frac13,\frac23]},\SP_n(\bar U,\b\bar U))
\endCD
$$
in which the maps $r,r^\prime$ are the natural \ll restriction\rr maps, 
and the 
maps $s_1,s_2$ are given by scanning. (Here, $\Map$ indicates
continuous maps.)
The map of the theorem is closely related to the map on fibres
of the horizontal maps. To obtain the map of the theorem
one must modify the method as explained in \cite{Se2}, but we ignore
these modifications here as our purpose is merely to explain the main
points of the argument.
After this modification, the map $r$ becomes a 
quasifibration.  It is an elementary fact that the map 
$r^\prime$ is a fibration. The theorem will follow if we  prove that 
the maps $s_1,s_2$ of total spaces and of base spaces are 
homotopy  equivalences. 

In the case of $s_1$, this is trivially so.  To deal with $s_2$, we 
shall need 
some new notation. Let $R_x=\{z\in\C\ \vert\ \text{Re}\,z\ge 
x\}$, and 
$L_x=\{z\in\C\ \vert\ \text{Re}\,z\le x\}$.  
Then consider the commutative diagrams
$$
\CD
\SP_n(D_{[0,1]},D_{[0,\frac13]}\cup D_{[\frac23,1]})  @>>> 
\SP_n(D_{[0,1]},R_{\frac16}\cup D_{[0,\frac13]}\cup 
D_{[\frac23,1]})\\
@VVV @VVV\\
\SP_n(D_{[0,1]},L_{-\frac{1}6}\cup D_{[0,\frac13]}\cup 
D_{[\frac23,1]}) @>>>
\SP_n(D_{[0,1]},R_{\frac16}\cup L_{-\frac{1}6}\cup 
D_{[0,\frac13]}\cup 
D_{[\frac23,1]})
\endCD
$$
and
$$
\CD
\Map(D_{[\frac13,\frac23]},\SP_n(\bar U,\b\bar U)) @>>>
\Map(L_{\frac16}\cap D_{[\frac13,\frac23]},\SP_n(\bar U,\b\bar U))\\
@VVV @VVV\\
\Map(R_{-\frac{1}6}\cap D_{[\frac13,\frac23]},\SP_n(\bar U,\b\bar U)) 
@>>>
\Map(L_{\frac16}\cap R_{-\frac{1}6}\cap 
D_{[\frac13,\frac23]},\SP_n(\bar U,\b\bar U))
\endCD
$$
Both diagrams are homotopy Cartesian. There are compatible maps 
from the 
first diagram to the second, given by scanning, one of which is the 
map 
$s_2$. The remaining three maps are homotopy equivalences 
(because $s_1$ 
is). Hence $s_2$ must be a homotopy equivalence also, as required.
\qed
\enddemo

\noindent{\it The scanning construction for algebraic maps}

It is illuminating to convert Theorem 2.2 to a result about polynomial
functions with singularities of a certain type.  This gives
the connection with the work of Vassiliev (\cite{Va1}, \cite{Va2}) 
mentioned in the introduction.  We need the following definition:

\proclaim{Definition 2.3} (1) For $n\ge 2$, $\CalSP^d_n(\C)$
denotes the space of (not necessarily monic) complex 
polynomial functions $f(z)=\sum a_i z^i$ of degree exactly $d$,
such that every root of $f$ has multiplicity less than $n$. 

\noindent (2) For $n\ge 2$, and any nonempty open subset
$X\sub\C$, $\CalSP_n(X)$ denotes 
the space of complex polynomial functions $f(z)=\sum a_i z^i$ 
such that every root of $f$ in $X$ has multiplicity less than $n$,
and such that $f$ is not identically zero.
\endproclaim

\no Both $\CalSP^d_n(\C)$ and $\CalSP_n(X)$ are topologized as subspaces
of the space of all complex polynomials. Note that $\CalSP_n(\C)$ is
bijectively equivalent to the disjoint union
$\cup_{d\ge 0} \CalSP^d_n(\C)$, but these spaces are not homeomorphic
because $\CalSP_n(\C)$ is connected ---  roots of polynomials in 
$\CalSP_n(\C)$ are allowed to move to infinity and \ll disappear\rrr.

There is a version of the scanning map for $\CalSP^d_n(\C)$,
namely
$$\spreadlines{2\jot}
\scan:\CalSP^d_n(\C)\to \Map(\C,\CalSP_n(U)),
\quad  f\mapsto (z\mapsto f\vert_{U(z)})
$$
(where, as in the definition of the earlier scanning map, we use the
canonical identification $U(z)\cong U$).
Here, $\Map(A,B)$ denotes the space of continuous maps
from $A$ to $B$.  
The definition of $\CalSP^d_n(\C)$ suggests that we
consider as well the jet map:
$$\spreadlines{2\jot}
\jet:\CalSP^d_n(\C)\to \Map(\C,\C^n-\{0\}),
\quad f\mapsto(f,f^\prime,\dots,f^{(n-1)}).
$$
We shall describe a relation between the scanning map for
configurations and these two natural maps.  

The maps
$$\spreadlines{2\jot}
\align
p&:\CalSP^d_n(\C)\to \SP^d_n(\C)\\
q&:\CalSP_n(U)\to \SP_n(\bar U,\partial\bar U)
\endalign
$$
given by assigning to a polynomial function its
roots are the key to this relation.  
It is obvious that $p$ is a principal fibre
bundle, with fibre $\C^\ast$. Moreover, this is a
trivial bundle because there is a section, defined by assigning to a
point of $\SP^d_n(\C)$ the corresponding monic polynomial.
Similarly, the \ll fibre\rr of the map $q$ is the space of
all polynomials in $\CalSP_n(U)$ whose roots
lie outside $U$, and this is homotopy equivalent to $\C^\ast$.
We claim that $q$ is in fact a quasifibration. This may be proved
using the well known criterion of Dold and Thom, as in the proof
of a similar assertion in Lemma 3.3 of \cite{Se2}. Namely, we
filter the base space $\SP_n(\bar U,\partial\bar U)$ by the number
of points in $U$, and use the fact that $q$ is a (trivial)
fibre bundle over each successive difference in this filtration.
The Dold-Thom \ll attaching map\rr has the effect of multiplying polynomials
with no roots in $U$ by a fixed polynomial $z-\al$, where $\al$ lies
outside $U$. Since $\al$ may be moved continuously to $1$,
the corresponding map of $\C^\ast$ is a homotopy equivalence,
as required. 

The scanning maps for 
$\SP^d_n(\C)$ and $\CalSP^d_n(\C)$ are related by the
following commutative diagram. (Diagrams of this type will play
a central role in this paper.) 
$$
\CD
{\CalSP^d_n(\C)}    @>{\scan}>>   {\Map(\C,\CalSP_n(U))}   
@>{\jet_0}>>   {\Map(\C,\C^n-\{0\})} \\
@VVV   @VVV   @VVV \\
{\CalSP^d_n(\C)/\C^\ast}    @>>>   
{\Map(\C,\CalSP_n(U)/\C^\ast})   
@>>>   {\Map(\C,\C^n-\{0\}/\C^\ast )}\\
@V=VV   @V\simeq VV   @V\cong VV \\
{\SP^d_n(\C)}    @.   {\Map(\C,\SP_n(\bar U,\partial\bar U))} @.
{\Map(\C,\C P^{n-1})}
\endCD
$$
The first and second columns are induced by $p$ and $q$ respectively.  
The map \ll $\jet_0$\rr is induced by the map
$\CalSP_n(U)\to \C^n-\{0\}$,
$f\mapsto (f(0),f^\prime(0),\dots,f^{n-1}(0))$.
Note that the first row of the diagram is simply a factorization
of the jet map $\CalSP^d_n(\C)\to \Map(\C,\C^n-\{0\})$.
 
Taking into account the behaviour of the scanning map at
$\infty$, we see that the second row of this diagram gives a map
into $\Om^2_d\C P^{n-1}$, which we denote by
$$\spreadlines{2\jot}
j^d_n:\SP^d_n(\C)\to \Omega^2_d \C P^{n-1}.
$$
With all the necessary preparations behind us, we can now prove
our first main result concerning this map:

\proclaim{Theorem 2.4}  For $n\ge 3$, 
$$\spreadlines{2\jot}
\lim_{d\to\infty} j^d_n:
\lim_{d\to\infty}\SP^d_n(\C)\to \lim_{d\to\infty}\Omega^2_d \C P^{n-1}
\simeq \Omega^2_0 \C P^{n-1}
$$
is a homotopy equivalence. For $n=2$, $\lim_{d\to\infty} j^d_n$
is a homology equivalence.
\endproclaim

\demo{Proof} Consider the above commutative diagram, in which
$j^d_n$ is (essentially) the second row.  The first part of the
second row, i.e. the scanning map
$\SP^d_n(\C)\to\Map(\C,\CalSP_n(U)/\C^\ast)$,
gives rise to the map of Theorem 2.2. Hence it is
a homotopy equivalence (in the limit $d\to\infty$) if $n\ge 3$,
and a homology equivalence if $n=2$. 

We claim that the jet map
$$\spreadlines{2\jot}
\jet_0:\CalSP_n(U)\to \C^n-\{0\},\quad
f\mapsto(f(0),f^\prime(0),\dots,f^{(n-1)}(0))
$$ 
is a $\C^\ast$-equivariant weak homotopy equivalence. (This
implies that the second part of the second row is a 
weak homotopy equivalence, and hence that the map in the statement of the
theorem is a weak homotopy equivalence. But each spaces in this statement
has the homotopy type of a CW-complex, so the map is actually a homotopy
equivalence, and the proof of the theorem is complete.)
To prove the claim, we use the same direct argument as for Proposition 1
of \cite{Ha} to show that the inclusion 
$$\spreadlines{2\jot}
\CalSP_n(U)\to
\{ f(z)=\sum_{i\ge 0} a_i z^i \st
(f(0),f^\prime(0),\dots,f^{(n-1)}(0))\ne (0,\dots,0)\}.
$$
is a weak homotopy equivalence.
On replacing $a_i$ by $ta_i$ for $i\ge n$, and letting $t\to 0$,
we deduce that $\CalSP_n(U)$ is weakly homotopy equivalent to
$$\spreadlines{1\jot}
\{ f(z)=\sum_{i=0}^{n-1} a_i z^i \st
(a_0,a_1,2!a_2,\dots,(n-1)!a_{n-1})\ne (0,\dots,0)\}.
$$
The jet map is therefore equivalent to the map
$$\spreadlines{2\jot}
\C^n-\{0\}\to \C^n-\{0\},\quad
(a_0,\dots,a_{n-1})\mapsto (a_0,a_1,2!a_2,\dots,(n-1)!a_{n-1})
$$
which is certainly a weak homotopy equivalence. Moreover, 
all maps here are clearly $\C^\ast$-equivariant.
\qed
\enddemo

\no {\it Remark.} A consequence of (the last part of) this 
proof is that the space 
$\SP_n(\bar U,\partial\bar U))$, which appears in the scanning 
construction for configurations, is homotopy equivalent 
to $\C P^{n-1}$.

\noindent{\it Segal's theorem on rational functions}

We shall give a brief description of (the stable version of)
Segal's theorem (\cite{Se2}) on holomorphic maps, in the spirit of
the above discussion.

\proclaim{Definition 2.5}
For $n\ge 2$, let $\Q^{(n-1)}_d(\C)$ be the space of $n$-tuples
$(\al_1,\dots,\al_n)$ with $\al_i\in\SP^d(\C)$
such that $\al_1\cap\dots\cap\al_n=\emptyset$.
\endproclaim

\noindent Alternatively, $\Q^{(n-1)}_d(\C)$ is the space of $n$-tuples
$(p_1,\dots,p_n)$ of complex polynomials of degree $d$ which
are monic and coprime.  The polynomials $p_1,\dots,p_n$ may be
regarded as the homogeneous coordinates of a holomorphic
map $F=[p_1;\dots;p_n]$ from $\C P^1=S^2=\C\cup\infty$
to $\C P^{n-1}$. We have $F(\infty)=[1;\dots;1]$ and
$[F]=d\in\pi_2\C P^{n-1}$.  Conversely, it is well known that any
holomorphic map $F:\C P^1\to\C P^{n-1}$ satisfying the last
two conditions corresponds to an element of $\Q^{(n-1)}_d(\C)$.
This means that $\Q^{(n-1)}_d(\C)$ may be identified with
the space $\Hol_d(S^2,\C P^{n-1})$  of such maps.

There is a natural inclusion map
$i^d_n:\Q^{(n-1)}_d(\C) = \Hol_d(S^2,\C P^{n-1}) \to \Omega^2_d\C P^{n-1}$.
In \cite{Se2}, Segal proves:

\proclaim{Theorem 2.6}  For $n\ge 3$, 
$$\spreadlines{2\jot}
\lim_{d\to\infty} i^d_n:
\lim_{d\to\infty}\Q^{(n-1)}_d(\C)\to \Omega^2_0\C P^{n-1}
$$
is a homotopy equivalence. For $n=2$, $\lim_{d\to\infty} i^d_n$
is a homology equivalence.
\endproclaim

We sketch the proof here, making only a rearrangement of the
proof given in \cite{Se2}.

\proclaim{Definition 2.7} (1) For $n\ge 2$, $\Cal Q^{(n-1)}_d(\C)$
denotes the space of $n$-tuples $(p_1,\dots,p_n)$
of (not necessarily monic)
complex  polynomial functions of degree exactly $d$, such
that $p_1,\dots,p_n$ have no common root.

\noindent(2)  For $n\ge 2$, and any nonempty open subset
$X\sub\C$, $\Cal Q^{(n-1)}(X)$ denotes
the space of $n$-tuples $(p_1,\dots,p_n)$
of complex  polynomial functions such
that $p_1,\dots,p_n$ have no common root in $X$, and
such that no $p_i$ is identically zero.
\endproclaim

There are scanning maps 
$$\spreadlines{2\jot}
\Q^{(n-1)}_d(\C)\to
\Map(\C,\Q^{(n-1)}(\bar U,\partial \bar U))
$$ 
(in which the definition of $\Q^{(n-1)}(X,Y)$ is analogous to that
of $\SP_n(X,Y)$)
and 
$$\spreadlines{2\jot}
\Cal Q^{(n-1)}_d(\C)\to\Map(\C,\Cal Q^{(n-1)}(U)).
$$
The analogue of the jet map in the present situation is the inclusion
$\Cal Q^{(n-1)}_d(\C)\to \Map(\C,\C^n-\{0\})$. The
analogue of \ll $\jet_0$\rr is simply \ll evaluation at $0$\rrr.
The analogue of the earlier commutative diagram is:
$$
\CD
{\Cal Q^{(n-1)}_d(\C)}    @>{\scan}>>  
{\Map(\C,\Cal Q^{(n-1)}(U))}    @>{\text{eval}_0}>>   
{\Map(\C,\C^n-\{0\})}
\\ @VVV   @VVV   @VVV \\
{\Cal Q^{(n-1)}(\C)/(\C^\ast)^n}    @>>>   
{\Map(\C,\Cal Q^{(n-1)}(U)/(\C^\ast)^n})   
@>>>   {\Map(\C,\C^n-\{0\}//(\C^\ast )^n})\\
@V=VV   @V\simeq VV   @V\simeq VV\\
{\Q^{(n-1)}(\C)}    @.   {\Map(\C,\Q^{(n-1)}(\bar U,\partial\bar U))} @.
{\Map(\C,\C P^{n-1})} 
\endCD
$$ 
The only new feature here is that $(\C^\ast)^n$ does not act 
{\it freely} on $\C^n-\{0\}$; for this reason it is necessary to use the
homotopy quotient $\C^n-\{0\}//(\C^\ast )^n$ in order to ensure
that each vertical map is a fibration with fibre $(\C^\ast)^n$. (Note
however that $(\C^\ast)^n$ acts freely on $\Cal Q^{(n-1)}(U)$, because of
the condition that no component of any element of $\Cal Q^{(n-1)}(U)$ is
identically zero.)

The proof of Segal's theorem now proceeds as in the case of Theorem 2.4. 
First, the map under consideration is given by the second row of the
diagram (after imposing the basepoint condition at $\infty$).
Next, the scanning map 
$Q^{(n-1)}_d(\C)\to \Om^2_d Q^{(n-1)}(\bar U,\partial\bar U)$
is a homotopy equivalence (in the limit $d\to\infty$),
as in Theorem 2.2.  
And finally, we argue as in the proof of Theorem 2.4 that
eval${}_0$ gives a $(\C^\ast)^n$-equivariant homotopy equivalence
$\Cal Q^{(n-1)}(U)\simeq \C^n-\{0\}$. (There is an additional
complication in the present situation as we have excluded from 
$\Cal Q^{(n-1)}(U)$ those $(p_1,\dots,p_n)$ which have some
component identically zero. However, the removal of this infinite-codimensional
space of functions does not affect the homotopy type --- see section 4
of \cite{Se2}.)

\no{\it Remark.} It follows from the homotopy equivalence
$\Cal Q^{(n-1)}(U)\simeq \C^n-\{0\}$ that the space 
$\Q^{(n-1)}(\bar U,\partial\bar U)$ is homotopy equivalent
to the homotopy quotient $\C^n-\{0\}//(\C^\ast )^n$. This,
in turn, is homotopy equivalent to the \ll fat wedge\rr $W_n(\C P^{\infty})$
of $n$ copies of $\C P^{\infty}$ (see section 2 of \cite{Se2}).

The first row of the above diagram is simply the natural inclusion
map.  Thus, we have shown that the natural inclusion 
$\Hol_d(S^2,\C P^{n-1})\to \Map_d(S^2,\C P^{n-1})$
may be identified up to homotopy with the scanning map
$\Q^{(n-1)}_d(\C)\to\Om^2_d\Q^{(n-1)}(\bar U,\partial\bar U))$,
which is a homotopy equivalence (in the limit $d\to\infty$).
This completes our sketch of the proof of Segal's theorem.

Now we describe the relation between
Theorem 2.4 and Segal's theorem.
Theorem 2.4 says that a certain map
$$\spreadlines{2\jot}
j^d_n:\SP^d_n(\C)\to \Omega^2_d \C P^{n-1}
$$
is a homotopy equivalence (when $d\to\infty$).
Theorem 2.6 says that a certain map
$$\spreadlines{2\jot}
i^d_n:\Q^{(n-1)}_d(\C)\to \Omega^2_d \C P^{n-1}
$$
is a homotopy equivalence (when $d\to\infty$). 
We have a map
$$\spreadlines{2\jot}
T^d_n:\SP^d_n(\C)\to Q^{(n-1)}_d(\C),\quad
f\mapsto (f, f+f^\prime,\dots,f+f^{(n-1)})
$$
(this modification of the usual jet map ensures that the
right hand side is an $n$-tuple of monic polynomials
of degree exactly $d$).
All these maps are related by the following commutative diagram
$$
\CD
\SP^d_n(\C)  @>{j^d_n}>> \Omega^2_d\CalSP_n(U)\simeq\Om^2_d\C P^{n-1}\\
@V{T^d_n}VV  @VVV\\
\Q^{(n-1)}_d(\C) @>{i^d_n}>> 
\Omega^2_d\Cal Q^{(n-1)}_d(U)\simeq\Om^2_d\C P^{n-1}
\endCD
$$
where the right hand vertical map is induced by the (modified) jet map
$\CalSP_n(U)\to \Cal Q^{(n-1)}(U)$. Exactly as in the proof of Theorem 2.4,
one can show that this jet map is homotopic to the identity map.

It follows from this that we have a 
factorization $j^d_n\simeq i^d_n\circ T^d_n$. Hence:

\proclaim{Corollary 2.8} The limit (as $d\to\infty$) of
$T^d_n:\SP^n_d(\C)\to Q^{(n-1)}_d(\C)$
is a homotopy equivalence for $n\ge 3$, and a homology
equivalence for $n=2$.
\qed
\endproclaim

It turns out that there is a much more precise relation between
$\SP^n_d(\C)$ and $Q^{(n-1)}_d(\C)$. In the introduction 
we referred to Vassiliev's
result  (\cite{Va1}, \cite{Va2}) that 
$\SP^{nd}_n(\C)$ and $Q^{(n-1)}_d(\C)$ are stably homotopy
equivalent. To provide some motivation for this, 
we need \ll unstable\rr versions of
Theorems 2.4 and 2.6.  Regarding Theorem 2.4, we have:

\proclaim{Theorem 2.9} The map 
$$\spreadlines{2\jot}
j^d_n:\SP^d_n(\C)\to \Omega^2_d \C P^{n-1}
$$
(defined earlier)
is a homotopy equivalence up to dimension $(2n-3)[d/n]$
if $n\ge 3$, and 
a homology equivalence up to dimension $(2n-3)[d/n]$
if $n=2$.
\endproclaim

\demo{Proof} It follows from \cite{Ar}
that the stabilization map
$\SP^d_n(\C)\to \SP^{d+1}_n(\C)$
is a homology equivalence up to dimension $(2n-3)[d/n]$.
(Further details concerning this result are given in the Appendix.)
Taken together with Theorem 2.4, this gives the homology
statement of the theorem. If $n\ge 3$ both spaces are simply
connected, so we obtain the homotopy statement as well.
\qed
\enddemo

Regarding Theorem 2.6, Segal considered the stabilization map
$\Q^{(n-1)}_d(\C)\to \Q^{(n-1)}_{d+1}(\C)$ in section 5 of
\cite{Se2}.  As a consequence of this and Theorem 2.6, he
obtained the following unstable result:

\proclaim{Theorem 2.10} The inclusion map 
$$\spreadlines{2\jot}
i^d_n:\Q^{(n-1)}_d(\C)\to \Omega^2_d\C P^{n-1}
$$
is a homotopy equivalence up to dimension $(2n-3)d$ if $n\ge 3$,
and a homology equivalence up to dimension $(2n-3)d$
if $n=2$.
\qed
\endproclaim

The last two theorems show that $\SP^{nd}_n(\C)$ and
$\Q^{(n-1)}_d(\C)$ have isomorphic homology groups up to
dimension $(2n-3)d-1$, which indicates the plausibility of
Vassiliev's result.
Proofs of Vassiliev's stable homotopy equivalence between these
two spaces can be found in \cite{GKY2} and \cite{Kl2}.

$${}$$

\head
\S 3 Further examples of the same type
\endhead

\noindent{\it A theorem of Vassiliev}

We turn now to Vassiliev's original result (mentioned
in the introduction), which concerns {\it real} polynomials.

\proclaim{Definition 3.1} $\P^d_n(\R)$ denotes the space
of {\it real} polynomials of degree $d$ which are {monic},
and all of whose {\it real} roots have {multiplicity less than} $n$.
\endproclaim

\noindent The space $\P^d_n(\R)$ is a subspace of
$\SP^d(\C)$ (but it is not a subspace of $\SP^d_n(\C)$
because no conditions are imposed on the non-real roots).
We have a \ll horizontal scanning map\rr
$$\spreadlines{2\jot}
\P^d_n(\R)\to \Omega \P_n(I,\b I)
$$
where $I=[-1,1]$ and
$$\spreadlines{2\jot}
\P_n(I,\b I)=
\{\text{$\sum$} d_ix_i
\in \SP(I\times[-1,1],\b I\times[-1,1]\ \vert\ 
d_i<n\ \text{when}\ x_i\in\R\}.
$$
This is defined by associating to a configuration $\al$ the map
$$\spreadlines{2\jot}
x\mapsto \al\cap\bar V(x)\in \P_n(\bar V(x),\b \bar V(x))
\cong \P_n(I,\b I),
$$
where $V(x)=\{z\in\C \st \vert \Re(z)-x\vert < \ep \}$.

The analogue of Definition 2.3 here is:

\proclaim{Definition 3.2} (1) For $n\ge 2$, $\Cal P^d_n(\R)$
denotes the space of (not necessarily monic) real 
polynomial functions of degree exactly $d$, all of whose
real roots have multiplicity less than $n$.

\noindent (2) For $n\ge 2$, and any nonempty open subset
$X\sub\R$, $\Cal P^d_n(X)$
denotes the space of real polynomial functions $f(x)=\sum a_i x^i$ 
such that every root of $f$ in $X$ has multiplicity less than $n$,
and such that $f$ is not identically zero.
\endproclaim

\noindent As in the case of complex polynomials,
there is a horizontal scanning map
$$\spreadlines{2\jot}
\scan:\Cal P^d_n(\R)\to \Map(\R,\Cal P_n((-1,1)), \quad
f\mapsto (x\mapsto f\vert_{V(x)})
$$
and a jet map
$$\spreadlines{2\jot}
\jet:\Cal P^d_n(\R)\to \Map(\R,\R^n-\{0\}),
\quad  f\mapsto(f,f^\prime,\dots,f^{(n-1)})\vert_{\R}.
$$

We have a commutative diagram analogous to the 
ones in the last section:

$$
\CD
{\Cal P^d_n(\R)}    @>{\scan}>>   {\Map(\R,\Cal
P_n((-1,1)))}    @>{\jet_0}>>   {\Map(\R,\R^n-\{0\})} \\
@VVV   @VVV   @VVV \\
{\Cal P^d_n(\R)/\R^\ast}    @>>>   
{\Map(\R,\Cal P_n((-1,1))/\R^\ast})  
@>>>   {\Map(\R,\R^n-\{0\}/\R^\ast )}\\
@V=VV   @V\simeq VV   @V\simeq VV\\
{\P^d_n(\R)}    @.   {\Map(\R,\P_n([-1,1],\{-1,1\}))} @.
{\Map(\R,\R P^{n-1})} 
\endCD
$$  
Let $j^d_n:\P^d_n(\R)\to\Omega_d\R P^{n-1}$ be the
second row of this diagram. We can now give a proof of Vassiliev's result
concerning $\P^d_n(\R)$ as $d\to\infty$.

\proclaim{Theorem 3.3}
For $n\ge 4$, 
$$\spreadlines{2\jot}
\lim_{d\to\infty} j^d_n:
\lim_{d\to\infty}\P^d_n(\R)\to \Omega_0 \R P^{n-1}
$$
is a homotopy equivalence. For $n=3$, $\lim_{d\to\infty} j^d_n$
is a homology equivalence.
\endproclaim

\demo{Proof} This is exactly analogous to the proof of Theorem 2.4.
\qed
\enddemo

\noindent{\it Related examples}

There is a family of related examples which includes the spaces
$\SP^d_n(\Bbb C)$, $\P^d_n(\Bbb R)$:

\proclaim{Definition 3.4}  Let $X,Y$ be subspaces of $\C$.
Then $\P^d_{Y,n}(X)$ denotes the space of complex
monic polynomials $f$ of degree $d$ such that (i) $f(X)\sub X$
and (ii) $f$ has no $n$-fold roots in $Y$.
\endproclaim

\noindent We shall consider the four examples obtained by taking $X,Y=\R$
or $\C$.  
It will be convenient to express these results in terms of spheres
rather than projective spaces; note that $\Omega^2 S^{2n-1}$ may
be identified with the identity component of $\Omega^2 \C P^{n-1}$,
and that $\Omega S^{2n-1}$ may
be identified with the universal covering space of $\Omega \C P^{n-1}$.
First, we have
$\P^d_{\C,n}(\C)=\SP^d_n(\Bbb C)$ and $\P^d_{\R,n}(\R)=\P^d_n(\Bbb R)$,
and we have seen in Theorems 2.4 and 3.3 that these provide
finite-dimensional models for the loop spaces
$\Omega^2\Bbb CP^{n-1}$ and
$\Omega\Bbb RP^{n-1}$.

Concerning $\P^d_{\R,n}(\C)$ we have:

\proclaim{Theorem 3.5}
For $n\ge 3$, there is a homotopy equivalence
$$\spreadlines{2\jot}
\lim_{d\to\infty}
\P^d_{\R,n}(\C)\to \Omega S^{2n-1}.
$$
\endproclaim

\demo{Proof} The method of Theorem 3.3 (using a horizontal
scanning map)  applies equally well here.
\qed\enddemo

Concerning $\P^d_{\C,n}(\R)$, there is a similar statement,
but this time involving the involution $\th:\C\to \C$, $\th(z)=\bar z$
(and its natural extension to $\C^{n}$).

\proclaim{Theorem 3.6}
For $n\ge 3$, there is a homotopy equivalence
$$\spreadlines{2\jot}
\lim_{d\to\infty}
\P^d_{\C,n}(\R)\to (\Omega^2 S^{2n-1})^{\th},
$$
where $(\Omega^2 S^{2n-1})^{\th}$ is the subspace of
$\Omega^2 S^{2n-1}$
consisting of maps $f$ which satisfy the condition
$\th({f(z)})=f(\th({z}))$ 
for all $z\in \C\cup\infty$.
\endproclaim

\demo{Proof} The method of Theorem 2.4 can be modified by imposing
\ll $\th$-equivariance\rr at the appropriate
\footnote{
In fact, the validity of the equivariant scanning argument in the case of 
more general finite group actions has already been noted in \cite{Se3}.}
points.
\qed\enddemo

As we pointed out in the introduction, Theorems 2.4, 3.3, 3.5 and 3.6
may be summarized in a single statement. To avoid special cases which
occur for low values of $n$, we assume $n\ge 4$ in the following version.

\proclaim{Theorem 3.7}
There is a $\th$-equivariant homotopy equivalence 
$$
\lim_{d\to\infty}
\P^d_{Y,n}(\C) \to \Map(Y^+,S^{2n-1})
$$
(if $n\ge 4$), where 
$Y=\R$ or $\C$,  $Y^+=\R\cup\{\infty\}$ or $\C\cup\{\infty\}$,
and where $\Map$
indicates
the space of basepoint preserving continuous maps.
\qed
\endproclaim

There is an analogous family of results relating to Segal's space
$\Q^{(n-1)}_d(\C)$ (the definition of this space was given in \S 2).
We may introduce a space $\Q^{Y,(n-1)}_d(X)$ for any subsets $X,Y$ of $\C$,
in which the polynomials $p_1,\dots,p_n$  are required to satisfy
the modified conditions  (i) $p_i(X)\sub X$
and (ii) $p_1,\dots,p_n$ have no common root in $Y$. There are four
basic examples. First we have $\Q^{\C,(n-1)}_d(\C) = \Q^{(n-1)}_d(\C)$ and 
$\Q^{\R,(n-1)}_d(\R)$.  Then there are two 
\ll mixed\rr spaces, $\Q^{\C,(n-1)}_d(\R)$ and $\Q^{\R,(n-1)}_d(\C)$.
The first of these is the space of \ll real\rr algebraic maps
$\C P^1\to \C P^{n-1}$, and this space (for $n=2$ at least) has
already been discussed by Segal in \cite{Se2}. 

We have the following analogue of Theorem 3.7:

\proclaim{Theorem 3.8}
There is a $\th$-equivariant homotopy equivalence 
$$
\lim_{d\to\infty}
\Q_d^{Y,(n-1)}(\C) \to \Map(Y^+,S^{2n-1})
$$
(if $n\ge 4$), where 
$Y=\R$ or $\C$,  $Y^+=\R\cup\{\infty\}$ or $\C\cup\{\infty\}$, and where
$\Map$ indicates
the space of basepoint preserving continuous maps.
\qed
\endproclaim

\noindent{\it Spaces of holomorphic maps \ll with bounded multiplicities\rrr}

Segal's result on the space $\Hol(S^2,\C P^{n-1})$ may be generalized
by imposing conditions of bounded multiplicity on the various
polynomials involved.  The appropriate generalization of Definition
2.5 is:

\proclaim{Definition 3.9}
For $n,m\ge 2$, let $\Q^{(n-1),m}_d(\C)$ be the space of $n$-tuples
$(\al_1,\dots,\al_n)\in\SP^d_m(\C)\times\dots\times\SP^d_m(\C)$
such that $\al_1\cap\dots\cap\al_n=\emptyset$.
\endproclaim

\no We denote by $\Hol^m_d(S^2,\C P^{n-1})$ the subspace of
$\Hol_d(S^2,\C P^{n-1})$ which corresponds to $\Q^{(n-1),m}_d(\C)$.
Observe that we have a filtration
$$\spreadlines{2\jot}
\es=\Hol^1_d(S^2,\C P^{n-1})\sub\dots\sub
\Hol^{d+1}_{d}(S^2,\C P^{n-1})=\Hol_d(S^2,\C P^{n-1})
$$
Corresponding to Definition 2.7 we have:

\proclaim{Definition 3.10} (1) For $n,m\ge 2$, $\Cal Q^{(n-1),m}_d(\C)$
denotes the space of $n$-tuples $(p_1,\dots,p_n)$ with
$p_i\in \CalSP^d_m(\C)$ such that
$p_1,\dots,p_n$ have no common root.

\noindent(2)  For $n,m\ge 2$, and any nonempty open subset
$X\sub\C$, $\Cal Q^{(n-1),m}_d(X)$ denotes
the space of $n$-tuples $(p_1,\dots,p_n)$ with
$p_i\in \CalSP^d_m(X)$ such
that $p_1,\dots,p_n$ have no common root in $X$.
\endproclaim

We can now proceed in the usual way.  The only new task is to study
the jet map
$$\spreadlines{2\jot}\gather
\jet_0:\Cal Q^{(n-1),m}(U)\to A_{n,m}\\
(p_1,\dots,p_n)\mapsto 
\left(  (p_1(0),p_1^{\prime}(0),\dots,p_1^{(m-1)}(0)),
\dots, (p_n(0),p_n^{\prime}(0),\dots,p_n^{(m-1)}(0))  \right)
\endgather
$$
where
$$\spreadlines{2\jot}
A_{n,m}=\{(v_1,\dots,v_n)\in (\C^m-\{0\})^n
\st ((v_1)_1,\dots,(v_n)_1)\ne (0,\dots,0)\}.
$$
As in the proof of Theorem 2.4, it can be shown that $\jet_0$
is a $(\C^\ast)^n$-equivariant homotopy equivalence.  This leads to
the following generalization of Theorem 2.6, for the analogous
map $j^d_{n,m}:\Q^{(n-1),m}_d(\C)\to \Omega^2_d A_{n,m}$:

\proclaim{Theorem 3.11}  Let $n\ge 2$. For $m\ge 3$, 
$$\spreadlines{2\jot}
\lim_{d\to\infty} j^d_{n,m}:
\lim_{d\to\infty}\Q^{(n-1),m}_d(\C)\to \Omega^2_0 A_{n,m}
$$
is a homotopy equivalence. For $m=2$, $\lim_{d\to\infty} j^d_{n,m}$
is a homology equivalence.
\qed
\endproclaim

The subset $A_{n,m}/(\C^\ast)^n$ of $(\C P^{m-1})^n$ consists of
elements $([v_1],\dots,[v_n])$ such that at least one of 
$[v_1],\dots,[v_n]$ lies in the \ll big cell\rr
$\{[v]\in\C P^{m-1} \st (v)_1\ne 0 \} \cong \C^{m-1}$. Hence it is
homotopy equivalent to the fat wedge $W_n(\C P^{m-1})$, i.e.
the set of all $([v_1],\dots,[v_n])$ such that at least one of 
$[v_1],\dots,[v_n]$ is equal to the basepoint $[1;0;\dots;0]$
of $\C P^{m-1}$.  Thus, Theorem 3.11 gives a homotopy equivalence
$$
\lim_{d\to\infty}\Q^{(n-1),m}_d(\C)\to \Omega^2_0 W_n(\C P^{m-1}).
$$
This shows that the space $A_{n,m}$ has a natural interpretation,
which is consistent with Segal's original approach (for the case
$m=\infty$).

$${}$$

\head
\S 4 The h-Principle
\endhead

\noindent{\it Spaces of polynomials with 
roots of bounded multiplicity}

As preparation for the h-Principle, we shall give a
general formulation which includes all the examples considered so far in this
paper.  Let $\Cal F_D(\C)$ denote the space of $n$-tuples 
$(p_1,\dots,p_n)$ of (not identically zero) complex polynomials
satisfying conditions of the following types:

\no(i) \ll degree conditions\rr (e.g. that $\deg p_i=d_i$,
where $D=(d_1,\dots,d_n)$ is fixed)

\no(ii) \ll coprime conditions\rr (e.g. that certain subsets
$p_{i_1},\dots,p_{i_k}$ of $p_1,\dots,p_n$ have no common factor)

\no(iii)\ll bounded multiplicity conditions\rr (e.g. that
all roots of $p_i$ have multiplicity less than $m_i$, where
$M=(m_1,\dots,m_n)$ is fixed).

\no Concrete examples of spaces $\Cal F_D(\C)$ satisfying
conditions of type (i) and (ii) are provided
\footnote{More precisely, the space $\Cal F_D(\C)/(\C^\ast)^n$
of $n$-tuples of {\it monic} polynomials corresponds to
the space of {\it based} holomorphic maps.}
by the spaces
$\Hol_D(S^2,X_{\De})$, where $X_{\De}$ is a smooth toric variety
defined by a fan $\De$ (see \cite{GKY1}, \cite{Gu2}). These
examples may be modified in obvious ways by imposing
conditions of type (iii). We may also consider the analogues
of the earlier spaces  $\P^d_{Y,n}(X)$ or  $\Q^{Y,(n-1)}_d(X)$
for $\Cal F_D(\C)$.

We regard $(p_1,\dots,p_n)$ as a holomorphic map $\C\to \C^n-A$,
where $A$ is the union of linear subspaces corresponding to
condition (ii). In the case of $\Hol_D(S^2,X_{\De})$,
the arrangement $A$ depends only on the fan of $X$,
so we write $A=A_{\De}$.  We then have
$X_{\De}\cong \C^n-A_{\De} / (\C^\ast)^{n-\dim X_{\De}}$
for a certain (free) action of $(\C^\ast)^{n-\dim X_{\De}}$
(see \cite{Co} and \cite{Gu2}). 

For a nonempty subset $X\sub\C$, we define $\Cal F(X)$ in the
same way, except that conditions (ii) and (iii) now apply only
to the roots of $p_1,\dots,p_n$ which lie in $X$ (and condition
(i) is omitted).  Just as in \S 2
and \S 3, we have a jet map
$$\spreadlines{2\jot}
\jet:\Cal F_D(\C)\to \Map(\C,{\C}^{n^{\prime}}-A^{\prime}),
$$
where ${\C}^{n^{\prime}}-A^{\prime}$ is the \ll prolongation\rr of
${\C}^{n}-A$ determined by condition (iii).  This factors as the top
row of the following diagram:
$$
\CD
{\Cal F_D(\C)}    @>{\scan}>>   {\Map(\C,\Cal F(U))}   
@>{\jet_0}>>   {\Map(\C,{\C}^{n^{\prime}}-A^{\prime})} \\
@VVV   @VVV   @VVV \\
{\Cal F_D(\C)/(\C^\ast})^n    @>>>   
\Map(\C,\Cal F(U)/(\C^\ast)^n)   
@>>>   {\Map(\C,{\C}^{n^{\prime}}-A^{\prime}//(\C^\ast)^n )}
\endCD
$$  
The action of $(\C^\ast)^n$ on $\Cal F(U)$ is always free. If the
action of $(\C^\ast)^n$ on ${\C}^{n^{\prime}}-A^{\prime}$ is free, we
may replace the homotopy quotient ${\C}^{n^{\prime}}-A^{\prime}//(\C^\ast)^n$
by the ordinary quotient.

As in the proof of Theorem 2.4, we can show:

\proclaim{Lemma 4.1}  The map 
$\jet_0:\Cal F(U)\to {\C}^{n^{\prime}}-A^{\prime}$ is a
$(\C^\ast)^n$-equivariant homotopy equivalence.
\qed
\endproclaim

\no Because of the interpretation of $\Cal F(\C)/(\C^\ast)^n$ and
$\Cal F(U)/(\C^\ast)^n$ as configuration spaces, the scanning method
(as in Theorem 2.2) leads to:

\proclaim{Lemma 4.2}  The map 
$\scan/(\C^\ast)^n:\Cal F_D(\C)/(\C^\ast)^n\to 
\Om^2_D\,\Cal F(U)/(\C^\ast)^n$ is a
homotopy equivalence in the limit $D\to\infty$.
\qed
\endproclaim

We deduce (from the lemmas and the diagram) the following
generalization of Theorems 2.2 and 2.6:

\proclaim{Theorem 4.3} The jet map induces a map
$$\spreadlines{2\jot}
j_D:\Cal F_D(\C)/(\C^\ast)^n\to
\Om^2_D \, {\C}^{n^{\prime}}-A^{\prime},
$$
and $\lim_{D\to\infty}j_D$ is a homotopy equivalence
(or homology equivalence, if the roots of any polynomial
in the definition of $\Cal F_D(\C)$ are required by
condition (iii) to be distinct).
\qed
\endproclaim

\no In the special case where 
$\Cal F_D(\C)/(\C^\ast)^n=\Hol_D(S^2,X_{\De})$ for some toric
variety $X_{\De}$, we have ${\C}^{n^{\prime}}-A^{\prime}={\C}^{n}-A$
and $\Om^2_D \, {\C}^{n^{\prime}}-A^{\prime} \simeq 
\Om^2_D \, {\C}^{n}-A  \simeq \Om^2_D \, X_{\De}$. 
In this case the map $j_D$ of Theorem 4.3
may be identified with the natural inclusion
$\Hol_D(S^2,X_{\De}) \to \Om^2_D\,X_{\De}$.

In general, Theorem 4.3 exhibits $\Cal F_D(\C)/(\C^\ast)^n$
(which is a \ll space of polynomials with roots of 
bounded multiplicity\rrr) as a model for the double loop space of
${\C}^{n^{\prime}}-A^{\prime}$.  This, in turn, is homotopy
equivalent to the double loop space of a \ll generalized
wedge product\rr of finite-dimensional complex projective
spaces (see \cite{GKY1} and the remarks following Theorem 3.11
in the previous section).

\noindent{\it The h-Principle}

Let $M,N$ be smooth manifolds.  A smooth map $f:M\to N$ may
be regarded as a section of the trivial bundle $M\times N\to M$,
and its $k$-jet $j_k(f):M\to J_k(M,N)$ is then a section
of the $k$-jet bundle $J_k(M,N)\to M$. Thus we have a map
$$\spreadlines{2\jot}
j_k:\Map(M,N)\to \Sec(J_k(M,N)),\quad f\mapsto j_k(f)
$$
where $\Map$ and $\Sec$ denote smooth maps and smooth sections,
respectively.  This map is not in general surjective; an element
$s$ of $\Sec(J_k(M,N))$ is said to be {\it holonomic} or {\it integrable}
if there exists an element $f$ of $\Map(M,N)$ such that $j_k(f)=s$.

More generally, if $S$ is a closed subspace of $J_k(M,N)$, we
define
$$\spreadlines{2\jot}
\align
\Map^S(M,N)&=
\{f\in\Map(M,N)\st j_k(f)(M)\cap S=\es\}\\
\Sec^S(J_k(M,N))&=
\{s\in \Sec(J_k(M,N)) \st s(M)\cap S=\es\}.
\endalign
$$
Then we have
$$\spreadlines{2\jot}
j^S_k:\Map^S(M,N)\to \Sec^S(J_k(M,N)),\quad f\mapsto j_k(f),
$$
and the image of this map consists of all sections which satisfy
the integrability condition..

The Smale-Hirsch Principle, or the (parametrized) 
h-Principle of Gromov, says
that, {\it under certain conditions on $M,N$ and $S$, 
the map $j^S_k$ is a homotopy equivalence.}  Under such
favourable conditions, the integrability condition is
therefore \ll irrelevant from the point of view of topology\rrr.
An extensive treatment of the h-Principle and its
generalizations can be found in \cite{Gr}.  Even to summarize
this work briefly here would not be feasible; we just mention that
examples of situations where the Principal holds are often
found when $M$ and $N$ are open manifolds, or when
$S$ is not too large.

For example, let $\dim M<\dim N$, and let $J_1(M,N)-S$ be
defined by the condition that the derivative has maximal rank.
Thus, $\Map^S(M,N)$ is the space of smooth immersions of
$M$ in $N$.  Smale and Hirsch studied immersions of
$M=S^m$ in $N=\R^n$ and discovered that regular homotopy
classes of such immersions are in one-to-one correspondence
with the elements of $\pi_m V_m(\R^n)$, where $V_m(\R^n)$
is the Stiefel manifold of $m$-frames in $\R^n$. This is
consistent with the h-Principle, as $\Sec^S(J_1(M,N))$ is
homotopy equivalent to $\Map(S^m,V_m(\R^n))$ in this case.

As another example, let $M$ and $N$ be complex manifolds, and
take $J_1(M,N)-S$ to be defined by the condition that the derivative
is $\C$-linear. In this case we have $\Map^S(M,N)=\Hol(M,N)$.
On the other hand the space $\Sec^S(J_1(M,N))$ is homotopy
equivalent to $\Map(M,N)$ (since the space of $\C$-linear
transformations is an $\R$-linear subspace of the space of
$\R$-linear transformations). So the h-Principle holds if and
only if the inclusion $\Hol(M,N)\to\Map(M,N)$ is a
homotopy equivalence. This is certainly false in general
(see \cite{Gr}), in particular when $M=N=S^2$.  However,
Segal's theorem indicates that something can be salvaged
in this case.
In the remainder of this section we shall describe how
Theorem 4.3 may be approached via the h-Principle.

Our starting point is Vassiliev's observation that
$$\spreadlines{2\jot}
\lim_{D\to\infty} \Cal F_D(\C)
$$
is (weakly) homotopy equivalent to the space of {\it smooth} maps
$f:\C\to \C^n-A$ such that the image of  $\jet(f)$  lies
in ${\C}^{n^{\prime}}-A^{\prime}$ and such that $\jet(f)$
satisfies the same \ll condition at $\infty$\rr as elements
of $\Cal F_D(\C)$. Let us denote this space by $\F^\ast(\C)$.
Theorem 4.3 is therefore equivalent to the statement that
the jet map
$$\spreadlines{2\jot}
j:\F^\ast(\C)/(\C^\ast)^n \to \Om^2_0\, {\C}^{n^{\prime}}-A^{\prime}
$$
is a homotopy equivalence.  This statement, which concerns only smooth maps,
may be deduced from the h-Principle, as we shall now explaim.

Let us denote by $F(\C)$ (and similarly for $F(X)$) the space 
of smooth maps
$f:\C\to \C^n-A$ such that the image of  $\jet(f)$  lies
in ${\C}^{n^{\prime}}-A^{\prime}$. {\it We claim that the h-Principle is
valid in this case,} at least for any surface $X$ which is constructed
by successively attaching two-dimensional 
disks to the two-dimensional disk $U$.
This may be proved by induction. For the case
$X=U$, the validity of the h-Principle is proved by \ll shrinking $U$ down to a
point\rrr. The inductive step is achieved by using the fact that the functor
$F$ converts (certain) cofibrations into fibrations 
(see the last part of our explanation of the
proof of Theorem 2.2, for a special case). The details of this
argument are explained, in greater generality, in \cite{Ha} and \cite{Po}.

Now we consider the usual commutative diagram, but this time for
the functor $F$ (rather than for algebraic maps):
$$
\CD
{\F^\ast(\C)}    @>{\scan}>>   {\Map^\ast(\C,\F(U))}   
@>{\jet_0}>>   {\Map^\ast(\C,{\C}^{n^{\prime}}-A^{\prime})} \\
@VVV   @VVV   @VVV \\
{\F(\C)/(\C^\ast})^n    @>{\scan/(\C^\ast})^n>>   
\Map^\ast(\C,\F(U)/(\C^\ast)^n)   
@>{\jet_0}/(\C^\ast)^n>>   
{\Map^\ast(\C,{\C}^{n^{\prime}}-A^{\prime}//(\C^\ast)^n )}
\endCD
$$  
The map $j$ is given by the bottom row. It suffices, therefore, to
show that the top row is a homotopy equivalence. (The conditions at $\infty$
in the top row are defined to be those which descend to the given conditions
at $\infty$ in the bottom row; the symbol $\Map^\ast$ indicates
that these conditions are in force.) The fact that the top row is a homotopy
equivalence requires two observations:

\no(1) the jet map $\F(U)\to {\C}^{n^{\prime}}-A^{\prime}$ is a
homotopy equivalence, and

\no(2) the scanning map $\F^\ast(\C)\to \Map^\ast(\C,\F(U))$ is a
homotopy equivalence.

We have already seen that (1) is true, as it begins the inductive
argument for the proof of the h-Principle for $\F$. To prove (2), we
consider the following commutative diagram
$$
\CD
\F(\{\vert z\vert < 2\})  @>>>  \F(\{1 <\vert z\vert < 2\}) \\
@VVV @VVV\\
\Map^\ast( \{\vert z\vert < 2 \},\F(U))  @>>>
\Map^\ast( \{ 1<\vert z\vert < 2\},\F(U))
\endCD
$$
in which the the horizontal maps are given by restriction and
the vertical maps are given by scanning. The horizontal maps
are fibrations, and --- by the h-Principle for the
cases where $X$ is a disk or an annulus --- the vertical maps
are homotopy equivalences. Hence the map of fibres is also
a homotopy equivalence. But this is the scanning map
$\F^\ast(U)\to \Map^\ast(U,\F(U))$, which is homotopic to the
map of (2).

$${}$$

\head
Appendix
\endhead

For completeness, we shall give here the proof of Arnold's result
on the homology of $\SP^d_n(\C)$, which was used in the proof of
Theorem 2.9.  Rather than quote directly from \cite{Ar} (where a
minor error occurs in the statement), we shall sketch an argument
along the lines of \cite{Va1}.

\proclaim{Theorem A1 (\cite{Ar}) }
For $n,d\geq 2$, 
the stabilization map 
$\SP^d_n(\C)\to \SP^{d+1}_n(\C)$
is a homology equivalence up to dimension
$N(d,n)$, where 
$$\spreadlines{2\jot}
N(d,n)=\cases
(2n-3)[d/n] & \text{ if }[d/n]<[(d+1)/n]
\\
\infty & \text{ if }[d/n]=[(d+1)/n]
\endcases
$$
\endproclaim

\demo{Proof}
For brevity we shall write
$\text{SP}^d_n=\text{SP}^d_n(\Bbb C)$ and
$\SP^d=\SP^d(\C)$, and we omit explicit mention of coefficients
in (co)homology.
Let $\Sigma^d_n\sub \text{SP}^d$ denote the discriminant variety
consisting of all polynomials $f\in \text{SP}^d$
which have at least one root of multiplicity $n$.
Since $\text{SP}^d_n=\text{SP}^d-\Sigma^d_n$,
Alexander duality gives
$$\spreadlines{2\jot}
H^k(\text{SP}^d_n)\cong \overline{H}_{2d-k-1}(\Sigma^d_n)
\quad
\text{ if }0<k<2d,
\tag $\ast$
$$
where we use the notation $\overline{H}_*(X)=H_*(\overline{X})$,
with $\overline{X}=X\cup\{\infty\}$
the one-point compactification of (a locally compact space)
$X$.

Let 
$I:\Bbb C \to \Bbb C^{[d/n]}$
be the Veronese embedding,
$I(z)=(z,z^2,z^3,\dots,z^{[d/n]})$.
Let $f\in \Sigma^d_n$. Assume that $f$
has at least $s$ distinct roots $z_1,z_2,\dots,z_s$ of
multiplicity $n$.
In this case, 
we denote by $\Delta (f;\{z_1,z_2,\dots, z_s\})$
the $(s-1)$-dimensional
$open$ simplex in $\Bbb C^{[d/n]}$ with vertices
$I(z_1),I(z_2),\dots,I(z_s)$.
(Note that since $s\leq [d/n]$, the points
$I(z_1),I(z_2),\dots,I(z_s)$ are in general
position.)

Define the geometric resolution $G=G(\Sigma^d_n)$ of $\Sigma^d_n$
by
$$\spreadlines{2\jot}
G=G(\Sigma^d_n)=
\bigcup_{f\in \Sigma_n^d,\{z_1,\dots,z_s\}}\{f\}\times
\Delta(f;\{z_1,z_2,\dots,z_s\})
\sub \Sigma^d_n\times \Bbb C^{[d/n]}
$$
Projection onto the first factor is a surjective open proper map
$G\to \Sigma^d_n$, and this extends naturally to a map
$\pi:\overline{G}\to \overline{\Sigma^d_n}$.
It is known that $\pi$ is a homotopy equivalence
([Va1]).

Define the subspaces $\{F_p\}_{p\geq 0}$ of $\overline{G}$ by
$$\spreadlines{2\jot}
F_p=
\{\infty\}\cup
(\bigcup_{f\in \Sigma_n^d,s\leq p}\{f\}\times 
\Delta(f;\{z_1,z_2,\dots,z_s\})).
$$
There is an increasing filtration
$$\spreadlines{2\jot}
F_0\sub F_1 \sub \dots \sub F_{[d/n]}=F_{[d/n]+1}=\dots
=\overline{G}\simeq
\overline{\Sigma^d_n}
$$
and so we have a homology spectral sequence 
$E_{p,q}^1=\overline{H}_{p+q}(F_p-F_{p-1})\Rightarrow
\overline{H}_{p+q}(\Sigma^d_n).$
 
If we take
$E^r_{p,q}=E_r^{p,2d-1-q}$, we obtain from $(\ast)$
a cohomology spectral sequence
$$\spreadlines{2\jot}
\Cal F^d_n=\{E_r^{p,q},d_r:E_r^{p,q}\to E_r^{p-r,q+1-r}\};
\quad
E_r^{p,q}\Rightarrow H^{q-p}(\text{SP}^d_n).
$$
Since there is a fibre bundle $F_p-F_{p-1}\to C_p(\Bbb C)$
with fibre homeomorphic to 
$\Bbb R^{2d-1-(2n-1)p}$ if $1\leq p\leq [d/n]$,
it follows from the Thom isomorphism theorem and 
Poincar\acuteaccent e  duality
that
$$\spreadlines{2\jot}
E_1^{p,q}=\cases
H^{(2-2n)p+q}(C_p(\Bbb C)) & \text{ if }1\leq p\leq [d/n]
\\
0 & \text{ otherwise.}
\endcases
$$
Similarly we have a cohomology spectral sequence
$$\spreadlines{2\jot}
\Cal F^{d+1}_n=
\{\text{ }^{\prime}E_r^{p,q},d_r^{\prime}:
\text{ }^{\prime}E_r^{p,q}\to \text{ }^{\prime}E_r^{p-r,q+1-r}\};
\quad
\text{ }^{\prime}E_r^{p,q}\Rightarrow H^{q-p}(\text{SP}^{d+1}_n)
$$
such that
$$\spreadlines{2\jot}
\text{ }^{\prime}E_1^{p,q}=\cases
H^{(2-2n)p+q}(C_p(\Bbb C)) & \text{ if }1\leq p\leq [(d+1)/n]
\\
0 & \text{ otherwise.}
\endcases
$$

Note that the stabilization map
$\text{SP}^d_n\to \text{SP}^{d+1}_n$
extends naturally to a map
$\text{SP}^d\to \text{SP}^{d+1}$
and this induces a map
$\Sigma^d_n\to \Sigma^{d+1}_n$.
This map
extends to the
open embedding $\Sigma^d_n\times \Bbb C \to \Sigma^{d+1}_n$
(up to homotopy),
which preserves the corresponding filtrations.
Because one-point compactification
is contravariant for open embeddings,
we obtain a map
$\overline{s}:\overline{\Sigma}^{d+1}_n
\to \overline{\Sigma^d_n\times \Bbb C}=
\overline{\Sigma}^d_n\wedge S^2$.
Hence
the stabilization map
$H^j(\text{SP}^{d+1}_n)\to H^j(\text{SP}^{d}_n)$
corresponds to the map
$$\spreadlines{2\jot}
\overline{H}_{2d+1-j}(\Sigma_n^{d+1})\to
\overline{H}_{2d-1-j}(\Sigma^d_n)
$$
which is the composition
$$\spreadlines{2\jot}
\overline{H}_{2d+1-j}(\Sigma^{d+1}_n)
@>\overline{s}_*>>
\overline{H}_{2d+1-j}(\Sigma^d_n\times \Bbb C)
@>\text{suspension}>\cong> \overline{H}_{2d-1-j}
(\Sigma^d_n).
$$
Since the above homomorphism preserves the filtrations, this
induces a homomorphism of spectral sequences
$\{\phi_r^{p,q}:\text{ }^{\prime}E_r^{p,q}\to
E_r^{p,q}:r\geq 1,(p,q)\in \Bbb Z \times \Bbb Z\}$.

Because the corresponding maps between filtrations
are natural, the diagram
$$
\CD
F_p^{\prime}-F_{p-1}^{\prime} @>\overline{s}>> F_p-F_{p-1}
\\
@VVV @VVV
\\
C_p(\Bbb C) @>=>> C_p(\Bbb C)
\endCD
$$
is commutative  (for $1\leq p\leq [d/n]$).
From the above description of $F_p-F_{p-1}$ as a bundle, we have:

(i) If $[d/n]=[(d+1)/n]$,
$\phi_1^{p,q}$ is an isomorphism for all $p,q$.

(ii)
If $[d/n]<[(d+1)/n]$,
$\phi_1^{p,q}$ is an isomorphism if $p\leq 0$, or if
$1\leq p\leq [d/n]$
and
$(2-2n)p+q\geq 0$.

We now apply the comparison theorem for spectral sequences.
For (i), the required result follows immediately,
so let us assume that $[d/n]<[(d+1)/n]$.
Then since
$$\spreadlines{2\jot}
(2-2n)p+q\geq 0\quad  \text{ if and only if }
\quad q\geq (2n-2)p,
$$
for $1\leq p\leq [d/n]$
the dimension
$q-p\geq (2n-3)p$ attains the maximal value
$(2n-3)[d/n]$.
Hence the induced homomorphism
$H^j(\text{SP}^{d+1})\to H^j(\text{SP}^d_n)$
is an isomorphism if $j\leq (2n-3)[d/n]$.
By the universal coefficients theorem, the same result is
valid for homology groups, so the proof of the theorem is complete.
\qed
\enddemo

\proclaim{Corollary A2 (\cite{Ar})}
If $[d/n]=[(d+1)/n]$ and $n\geq 3$, the stabilization map
$\text{SP}^d_n(\Bbb C)\to \text{SP}^{d+1}_n(\Bbb C)$
is a homotopy equivalence.
\qed
\endproclaim

\eightpoint
      
\noindent{\it Department of Mathematics, University
of Rochester, Rochester, New York 14627, USA
\newline
and 
\newline Department of Mathematics, Tokyo Metropolitan University,
Minami-Ohsawa 1-1, Hachioji-shi, Tokyo 192, Japan}      
                                             
\noindent{\it Department of Mathematics, 
Toyama International University,                             
Kaminikawa, Toyama 930, Japan}
                          
\noindent{\it Department of Mathematics, The
University of Electro-Communications, Chofu, 
Tokyo 182, Japan}

\no martin\@math.rochester.edu, martin\@math.metro-u.ac.jp
\newline
andrzej\@tuins.ac.jp
\newline
kohhei\@prime.e-one.uec.ac.jp


\newpage
\Refs   
     
\widestnumber\key{\bf CCMM}

\ref 
\key{\bf Ar}
\by V.I. Arnold
\paper Some topological invariants of algebraic functions 
\yr 1970 
\vol 21
\jour Trans. Moscow Math. Soc.
\pages 30--52
\endref

\ref
\key{\bf Bo}
\by C.-F. B\"odigheimer
\paper Stable splittings of mapping spaces
\inbook Algebraic Topology
Springer Lecture Notes in Math. 1286 
\eds H.R. Miller and D.C. Ravenel
\publ Springer
\yr 1987
\pages 174--187
\endref

\ref
\key{\bf CCMM}
\by F.R. Cohen, R.L. Cohen, B.M. Mann and R.J. Milgram
\paper The topology of rational functions and 
divisors of surfaces 
\vol 166 
\jour Acta Math.
\yr 1991 
\pages 163--221 
\endref

\ref\key{\bf CJS}\by R.L. Cohen, J.D.S. Jones, and G.B. Segal
\paper Floer's infinite dimensional Morse theory
and homotopy theory
\inbook The Floer Memorial Volume
\eds H. Hofer, C.H. Taubes, A. Weinstein, and E. Zehnder
\bookinfo Progress in Math. 133
\publ Birkh\"auser
\yr 1995
\pages 297--325
(also published in
Surikaisekikenkyusho Kokyuroku 883, Kyoto, 1994, 68--96)
\endref

\ref
\key{\bf Co}\by D.A. Cox
\paper The homogeneous coordinate ring of a toric variety
\jour J. Alg. Geom.
\vol 4
\yr 1995
\pages 17--50
\endref

\ref
\key {\bf Fu}
\by K. Fukaya
\paper Topological field theory and Morse theory
\jour Sugaku Expositions
\paperinfo (translation from Suugaku, 46 (1994), 289--307)
\vol 10
\yr 1997
\pages 19--39
\endref

\ref
\key{\bf GM}
\by M. Goresky and R. MacPherson
\book  Stratified Morse Theory
\yr 1988
\publ Springer
\endref

\ref
\key{\bf Gr}
\by M. Gromov
\book  Partial Differential Relations
\yr 1986
\publ Springer
\endref

\ref
\key{\bf GKY1}\by M.A. Guest, A. Kozlowski, and K. Yamaguchi
\paper The topology of spaces of coprime polynomials           
\jour Math. Z.
\vol 217 \yr 1994 \pages 435--446
\endref

\ref
\key{\bf GKY2}
\by M.A. Guest, A. Kozlowski, and K. Yamaguchi
\paper Stable splitting of the space of polynomials with 
roots of bounded multiplicity
\jour J. Math. Kyoto Univ
\paperinfo to appear
\endref

\ref
\key{\bf Gu1}
\by M.A. Guest
\paper On the space of holomorphic maps from the 
Riemann sphere to the quadric cone 
\jour Quarterly Jour. of Math. 
\vol 45 
\yr 1994
\pages 57--75
\endref

\ref
\key{\bf Gu2}
\by M.A. Guest
\paper  The topology of the space of rational 
curves on a toric variety
\jour Acta Math. 
\vol 174 
\yr 1995
\pages 119--145
\endref

\ref
\key{\bf Ha}
\by A. Haefliger
\paper Lectures on the Theorem of Gromov
\inbook Proc. of Liverpool Singularities Symposium II,
Springer Lecture Notes in Math. 209
\ed C.T.C. Wall
\publ Springer
\yr 1970
\pages 128--141
\endref

\ref
\key{\bf JS}
\by I.M. James and G.B. Segal
\paper On equivariant homotopy type
\jour Topology
\vol 17
\yr 1978
\pages  267--272
\endref

\ref
\key{\bf Kl1}
\by S. Kallel
\paper Particle spaces on manifolds and generalized 
Poincar\acuteaccent e dualities
\paperinfo preprint
\endref 

\ref
\key{\bf Kl2}
\by S. Kallel
\paper The topology of spaces of maps from a Riemann surface
into complex projective space
\paperinfo preprint
\endref

\ref
\key{\bf Kt}
\by F. Kato
\paper 
\paperinfo Master's Thesis, Shinshu University, 1994 (in Japanese)
\endref 

\ref\key {\bf Mc1}
\by D. McDuff
\pages 91--107
\paper Configuration spaces of positive and negative particles
\yr 1975\vol 14\jour Topology                                 
\endref 

\ref
\key {\bf Mc2}
\by D. McDuff
\pages 88--95
\paper Configuration spaces 
\inbook K-Theory and Operator Algebras, Springer Lecture Notes
in Math. 575
\publ Springer
\yr 1977                                
\endref

\ref
\key{\bf Po}
\by V. Poenaru
\paper Homotopy theory and differential singularities
\inbook Manifolds -- Amsterdam 1970,
Springer Lecture Notes in Math. 197
\ed N.H. Kuiper
\publ Springer
\yr 1970
\pages 106--133
\endref
                              
\ref\key{\bf Se1}\by G.B. Segal
\paper Configuration spaces and iterated loop spaces
\jour Invent. Math.\vol 21\yr 1973\pages 213--221   
\endref                                          

\ref\key {\bf Se2}\by G.B. Segal\pages 39--72
\paper The topology of spaces of rational functions
\yr 1979\vol 143                                   
\jour Acta Math.
\endref 

\ref\key {\bf Se3}\by G.B. Segal
\paper Some results in equivariant homotopy theory
\paperinfo unpublished manuscript
\endref 

\ref
\key{\bf Va1}
\by V.A. Vassiliev
\book Complements of Discriminants of Smooth Maps: Topology and Applications
\yr 1992 (revised edition 1994)
\vol 98
\publ Translations of Math. Monographs, Amer. Math. Soc. 
\endref

\ref
\key{\bf Va2}
\by V.A. Vassiliev
\paper Topology of discriminants and their complements
\inbook Proc. Int. Cong. Math. 1994
\yr 1995
\publ Birkh\"auser
\pages 209--226
\endref
  
\endRefs
\enddocument